\documentclass{birkmult}
\pdfoutput=1
\usepackage{amscd}
\usepackage{amssymb}
\usepackage{pdfpages}
\usepackage{graphicx}

\newcommand{\C}{\mathbb C}
\newcommand{\Q}{\mathbb Q}

\newcommand{\X}{\mathcal X}
\usepackage[curve]{xypic}
\newcommand{\comment}[1]{\marginpar{\sffamily{\noindent\tiny #1
   \par}\normalfont}}
\renewcommand{\comment}[1]{}
\newbox\mybox
\def\overtag#1#2#3{\setbox\mybox\hbox{$#1$}\hbox to
  0pt{\vbox to 0pt{\vglue-#3\vglue-\ht\mybox\hbox to \wd\mybox
      {\hss$\scriptstyle#2$\hss}\vss}\hss}\box\mybox}
\def\undertag#1#2#3{\setbox\mybox\hbox{$#1$}\hbox to 0pt{\vbox to
    0pt{\vglue#3\vglue\ht\mybox\hbox to \wd\mybox
      {\hss$\scriptstyle#2$\hss}\vss}\hss}\box\mybox}
\def\lefttag#1#2#3{\hbox to 0pt{\vbox to 0pt{\vss\hbox to
      0pt{\hss$\scriptstyle#2$\hskip#3}\vss}}#1}
\def\righttag#1#2#3{\hbox to 0pt{\vbox to 0pt{\vss\hbox to
      0pt{\hskip#3$\scriptstyle#2$\hss}\vss}}#1}

\def\Dot{\lower.2pc\hbox to 2.5pt{\hss$\bullet$\hss}}
\def\Circ{\lower.2pc\hbox to 2.5pt{\hss$\circ$\hss}}
\def\Vdots{\raise5pt\hbox{$\vdots$}}
\def\splicediag#1#2{\xymatrix@R=#1pt@C=#2pt@M=0pt@W=0pt@H=0pt}

\renewcommand\frame[2][3pt]{\hbox{$\vcenter{\hbox{\vrule\vbox 
{\hrule\kern#1\hbox{\kern#1$#2$\kern#1}\kern#1\hrule}\vrule}}$}}
\newcommand\lineto{\ar@{-}}
\newcommand\dashto{\ar@{--}}
\newcommand\dotto{\ar@{.}}
\newcommand{\bt}{\bullet}

\newtheorem*{BStheorem}{Bhupal-Stipsicz Theorem}

\newtheorem{theorem}{Theorem}[section]

\newtheorem*{theorem*}{Theorem}
\newtheorem{proposition}[theorem]{Proposition}
\newtheorem{lemma}[theorem]{Lemma}
\newtheorem{corollary}[theorem]{Corollary}
\newtheorem*{corollary*}{Corollary}

\newtheorem{conjecture}[theorem]{Conjecture}
\newtheorem*{conjecture*}{Conjecture}

\theoremstyle{definition}
\newtheorem{remark}[theorem]{Remark}
\newtheorem*{example*}{Example}
\newtheorem*{examples*}{Examples}
\newtheorem*{remarks*}{Remarks}
\newtheorem{definition}[theorem]{Definition}

\begin{document}
\title[Complex surface singularities with rational homology disk smoothings]
{Complex surface singularities with rational homology disk smoothings}
\author{Jonathan Wahl}
\dedicatory{To Andr\'as N\'emethi on his 60th birthday}
\address{Department of Mathematics\\The University of North
  Carolina\\Chapel Hill, NC 27599-3250} \email{jmwahl@email.unc.edu}
\begin{abstract} A cyclic quotient singularity of type $p^2/pq-1$ ($0<q<p, (p,q)=1$) has a smoothing whose Milnor fibre is a  $\Q$HD, or rational homology disk (i.e., the Milnor number is $0$) (\cite{w3}, 5.9.1).  In the 1980's, we discovered additional examples of such singularities: three triply-infinite and six singly-infinite families, all weighted homogeneous.  Later work of Stipsicz, Szab\'{o}, Bhupal, and the author (\cite {SSW}, \cite {bs}) proved that these were the only weighted homogeneous examples.  In his UNC PhD thesis (unpublished but available at \cite{jf}), our student Jacob Fowler completed the analytic classification of these singularities, and counted the number of smoothings in each case, except for types $\mathcal W$, $\mathcal N$, and $\mathcal M$.  In this paper, we describe his results, and settle these remaining cases; there is a unique $\Q$HD smoothing component except in the cases of an obvious symmetry of the resolution dual graph.  The method involves study of configurations of rational curves on projective rational surfaces. \end{abstract}
\maketitle
    
\section{Introduction}    
  
Let $(X,0)$ be the germ of a complex normal surface singularity.
A  \emph{smoothing of $(X,0)$} is a morphism $f:(\mathcal X,0)\rightarrow (\mathbb C,0)$, where $(\mathcal X,0)$ is an isolated Cohen-Macaulay singularity, equipped with an isomorphism $(f^{-1}(0),0)\simeq (X,0)$.  The Milnor fibre $M$ of a smoothing is a general fibre $f^{-1}(\delta)$, a four-manifold whose boundary is the link of $(X,0)$. The first betti number of $M$ is $0$ \cite{gst}.  We say $f$ is a  $\Q$HD (or \emph{rational homology disk}) smoothing if the second betti number of $M$ is $0$ as well (the \emph{Milnor number} $\mu = 0$).  In such a case, $(X,0)$ must be a rational singularity.
 
   The basic examples are smoothings of the cyclic quotient singularities of type $p^2/pq-1$, where $0<q<p, (p,q)=1$ (\cite{w2}, (2.7)).   For $f(x,y,z)=xz-y^p$, one has that $f:\C^3\rightarrow \mathbb C$ is a smoothing of the $A_{p-1}$ singularity, whose Milnor fibre $M$ has Euler characteristic $p$.   Now consider the cyclic subgroup $G\subset GL(3,\C)$ generated by the diagonal matrix  $[\zeta, \zeta^q, \zeta^{-1}]$, where $\zeta=e^{2\pi i/p}$. $G$ acts freely on $\C^3-\{0\}$ and $f$ is $G$-invariant; so there is  a map $f:\C^3/G\equiv \mathcal X\rightarrow \C$, a smoothing of the cyclic quotient singularity $A_{p-1}/G$, which has type $p^2/pq-1$. The new Milnor fibre is the free quotient $M/G$, of Euler characteristic $1$, hence Milnor number $0$.
    
    More examples can be produced by a similar ``quotient construction"(\cite{w3}, 5.9.2).  For instance, let $f(x,y,z)=xy^{p+1}+yz^{q+1}+zx^{r+1}$, $N=(p+1)(q+1)(r+1)+1$, and $G\subset GL(3)$ the diagonal cyclic subgroup generated by $[\zeta, \zeta^{(q+1)(r+1)},\zeta^{-(r+1)}]$, where $\zeta=e^{2\pi i/N}$.  The resulting class of examples was later named $\mathcal W(p,q,r)$ in \cite{SSW}, (8.3).  Another class $\mathcal N(p,q,r)$ was obtained by replacing $(\C^3,0)$ by a hypersurface $(V,0)\subset (\C^4,0)$,  $G$ by a group of automorphisms acting freely off the origin, and an appropriate $f$.  Three more families are constructed in \cite{w5}.

       However, the major way to produce examples uses H. Pinkham's general method of ``smoothing with negative weight"  \cite{p3} for a weighted homogeneous singularity $(X,0)$.  Writing $X=$Spec $A$, where $A$ is a graded ring, form the $\C^*$-compactification  $\bar{X}=$Proj $A[t]$ (where $t$ has weight $1$).  $\bar{X}$ has a smooth curve $\bar{C}$=Proj $A$ at infinity (which we assume is rational), along which are several cyclic quotient singularities.  Resolving those singularities yields $\bar{X}'$ with a star-shaped collection of curves $\bar{E}$, consisting of $\bar{C}$ plus chains of rational curves.  The associated graph $\Gamma$ of these curves is ``dual'' to the star-shaped resolution graph $\Gamma'$ of the singularity; it is non-degenerate, of signature $(1,s)$. (See e.g.\cite{SSW},(8.1) for details.)  (Beware: $\Gamma'$ is itself sometimes called the ``dual resolution graph.'')  A \emph{smoothing of negative weight of $X$} is a smoothing which can be extended to a smoothing of $\bar{X}'$ to which $\bar{E}$ lifts and is deformed trivially.  The general fibre is a smooth projective surface $Z$, with $H^1(\mathcal O_Z)=0$, containing a curve $E$ isomorphic to $\bar{E}$, which supports an ample divisor.  The Milnor fibre $M$ of this smoothing may be identified with the affine variety $Z-E$ (\cite{w3}, (2.2)).  Thus, $M$ is a $\Q$HD if and only if the curves of $E$ are rational and form a rational basis of Pic $Z$.  
    
Conversely, Pinkham shows how to construct a smoothing of negative weight of $(X,0)$ by starting with certain surface pairs $(Z,E)$ satisfying some cohomological vanishing. The author used this method to compile a large list of (only partially published) examples of $\Gamma$ which led to $\Q$HD smoothings.   
      The paper \cite{SSW} limited greatly the possible resolution graphs $\Gamma'$ of \emph{any} singularity admitting a $\Q$HD smoothing, and gave names to the author's families of examples (modified slightly in \cite{jf}). This work culminated in the Bhupal-Stipsicz theorem \cite{bs}, showing that the author's list of resolution dual graphs was complete for the weighted homogeneous case.         
       \begin{BStheorem} The resolution graphs (or dual graphs) of weighted homogeneous surface singularities admitting a $\Q$HD smoothing are exactly those of the following types: $p^2/(pq-1)$ cyclic quotients;   $\mathcal W(p,q,r)$; $\mathcal N(p,q,r)$; $ \mathcal M(p,q,r)$; $ \mathcal B^3_2(p$); $\mathcal C^3_2(p)$; $\mathcal C^3_3(p)$; $\mathcal A^4(p)$; $\mathcal B^4(p)$; $\mathcal C^4(p).$       
 \end{BStheorem}
 
  Resolution graphs and dual graphs for these singularities are listed at the end of this paper in Tables $A.1$ and $A.2$,  from the thesis of Jacob Fowler \cite{jf}. A node (or bullet) with no decoration is always assumed to be a $-2$ curve.  For the remainder of the paper, we disregard the well-understood cyclic quotients.
  
         Previous work by H. Laufer \cite{l} shows that the examples above with a central curve of valency $3$ are taut, i.e., have a unique analytic type, necessarily weighted homogeneous; further, all deformations (in particular, smoothings) are of negative weight.  In case the valency is $4$ (the last $3$ families), \cite{l} implies that the only analytic invariant is the cross-ratio of the central curve, and again all deformations are of negative weight.   
        
          Fowler's Ph.D. Thesis \cite{jf} attacked the key questions remaining for these $\Q$HD smoothings:
 \begin{itemize}
 \item Show the cross-ratios in the three infinite families of valency $4$ examples are uniquely determined, as in \cite{w5}.
 \item Determine the number of $\Q$HD smoothing components in each case.
 \item Calculate the fundamental groups of the Milnor fibres.
 \end{itemize}             
   
 We fix some language and notation.  Let $(X,0)$ be a weighted homogeneous surface singularity, of resolution dual graph $\Gamma$, admitting $\Q$HD smoothings.
    \begin{definition}  A \emph{$\Gamma$ surface} is a pair $(Z,E)$ consisting of a smooth rational surface and rational curve configuration $E$ such that
    the classes of the components of $E$ form a rational basis of $\text{Pic}\  Z$, and one is given an identification of the curve configuration $E$ with the graph $\Gamma$.
    \end{definition}
 \bigskip   
    
  Pinkham's Theorem in the current situation may be found in \cite{SSW}, (8.1)  and \cite{jf}, (2.2.3), yielding          
       \begin{theorem} (\cite{jf}, (2.3.1) ) Let $(X,0)$ be a singularity as above with a $\Q$HD smoothing, and resolution dual graph $\Gamma$.  Then there exists a one-to-one correspondence between $\Q$HD smoothing components of $(X,0)$ and $\Gamma$ surfaces $(Z,E)$ up to isomorphism.
       \end{theorem}   
       
Examples of $\Gamma$ surfaces are made as follows: take a specific curve configuration $D\subset \mathbb P^2$, blow up several times, obtaining $\pi:Z\rightarrow \mathbb P^2$, with $\pi^{-1}(D)$ consisting of a curve $E$ of type $\Gamma$ plus some $-1$ curves.  If $E$ spans Pic($Z$) rationally, one has a $\Gamma$ surface.  Given the location of the $-1$ curves in relation to the components of $E$, one can reverse the process and blow back down.

    For each $\Gamma$, Fowler makes very judicious choices of $D$ and the points to blow up, resulting in either one or two \emph{Basic Models} $(Z,E)$.  The models for $\mathcal W, \mathcal N$, and $ \mathcal M$ are listed on the first page of Table $A.2$, where small circles and light lines indicate the location of three $-1$ curves which allow the entire graph to be blown down, to four lines in general position.  Basic Models for the other graphs are more complicated and found in Fowler's thesis \cite{jf}.  The curve configuration $D$ will be unique up to projective equivalence.  One may get two models for the same $\Gamma$ and $D$ by blowing up in different ways (sometimes complex conjugate points).  The goal is to prove 
       
  \begin{conjecture}  Every $\Gamma$ surface $(Z,E)$ is a Basic Model.
  \end{conjecture}      
  
  In \cite{jf}, Fowler proves most of this Conjecture; his nearly complete result, explained below in Section $6$, states:
  
  \begin{theorem}\cite{jf}.  Suppose a $\Gamma$ surface $(Z,E)$ has self-isotropic subgroup which is basic.  Then $(Z,E)$ is a Basic Model.
  \end{theorem}
  Fowler also proves that the ``basic self-isotropic subgroup'' condition is automatically satisfied in all cases except for some $\Gamma$ of type $\mathcal W, \mathcal N$, or $ \mathcal M.$ 
  \begin{corollary}\cite{jf}  For $\Gamma$ not of type $\mathcal W, \mathcal N, \mathcal M$, every $\Gamma$ surface is a Basic Model.  In particular, for each valency $4$ example, there is a unique cross-ratio for which the corresponding singularity has a $\Q$HD smoothing.
  \end{corollary}
  
  The new contribution of the current paper is to handle the remaining cases.
  
  \begin{theorem}  Every $\Gamma$ surface of type $\mathcal W, \mathcal N, \mathcal M$ is a Basic Model.
  \end{theorem}
  
  The Basic Models for types $\mathcal W, \mathcal N$, and $\mathcal M$ start with four lines in general position, for which the fundamental group of the complement is abelian.  Therefore the Milnor fibre of a $\Q$HD smoothing of a singularity of this type has abelian fundamental group (hence is easily computable from $\Gamma$).  More generally, we can conclude
  
  \begin{theorem} Let $M$ be the Milnor fibre of a $\Q$HD smoothing of a singularity of type $\Gamma$.
  \begin{enumerate}
  \item If $\Gamma$ is of type $\mathcal W, \mathcal N$, or $\mathcal M$, then $\pi_1(M)$ is abelian.
  \item If $\Gamma$ is of type $\mathcal A^4$, $\mathcal B^4$, or $\mathcal C^4,$ then $\pi_1(M)$ is metacyclic, as described in \cite{w5}.
  \end{enumerate}
  \end{theorem}
  
  It remains an open problem to determine whether the fundamental group is abelian in cases $ \mathcal B^3_2, \mathcal C^3_2, \mathcal C^3_3$.
  
  Once we know the number of $\Q$HD smoothings from the main theorem, we can conclude that the explicit examples from the quotient construction (for types $\mathcal W, \mathcal N$, $\mathcal A^4$, $\mathcal B^4, \mathcal C^4$) give a complete list of smoothing components in those cases.  This also gives the only way to compute the metacyclic fundamental group.
  
   In Section $6$, we give more details on Fowler's method and list the number of $\Q$HD smoothing components for each $\Gamma$. 
   
   The isomorphism type of a $\Gamma$ surface of type $\mathcal W, \mathcal N$, or $ \mathcal M$  is determined by the location of the three extra $-1$ curves that are attached to $E$.  For special values of $p,q,r$ for which the graph $\Gamma$ has a symmetry, there could be a second location of the $-1$'s, leading to a different $\Gamma$ surface (which, one recalls, comes equipped with a specific identification with the graph).  In light of our new result, we find
    \begin{theorem} \cite{jf} Consider $\Q$HD smoothing components for type $\mathcal W, \mathcal N$, and $\mathcal M$.
  \begin{enumerate}
  \item There are two components for $\mathcal W(p,p,p)$, $\mathcal N(q+2, q,0)$, and $\mathcal M(r+1,q,r)$
  \item In all other cases, there is a unique $\Q$HD smoothing component.
   \end{enumerate}
   \end{theorem}    
  
    (Actually, \cite{jf} neglected to mention the exceptional $\mathcal N$ case, but it fits in easily with his work.)
     
 Our method, already used in \cite{bs}, is to blow up and down the given $\Gamma$ surface $(Z,E)$ so that one obtains a surface $(Z',E')$ with central curve of self-intersection $+1$, from which a blowing-down map to $\mathbb P^2$ can be constructed.  We analyze the possible singularities of the image of $E'$ and the blowing up needed to reach back to $Z'$, leading to location of all possible sets of essential $-1$ curves on $Z$ needed for blowing down.   All solutions will be Basic Models.

\section{Locations of $-1$ curves}

Suppose $\Gamma$ is a graph of smooth rational curves $E=\Sigma E_i$:

$$\xymatrix@R=4pt@C=24pt@M=0pt@W=0pt@H=0pt{\\
\lefttag{\bullet}{n_2/q_2}{8pt}\dashto[ddrr]&
&\hbox to 0pt{\hss\lower 4pt\hbox{.}.\,\raise3pt\hbox{.}\hss}
&\hbox to 0pt{\hss\raise15pt
\hbox{.}\,\,\raise15.7pt\hbox{.}\,\,\raise15pt\hbox{.}\hss}
&\hbox to 0pt{\hss\raise 3pt\hbox{.}\,.\lower4pt\hbox{.}\hss}
&&\righttag{\bullet}{n_{t-1}/q_{t-1}}{8pt}\dashto[ddll]\\
\\
&&\bullet\lineto[dr]&&\bullet\lineto[dl]\\
\lefttag{\bullet}{n_1/q_1}{8pt}\dashto[rr]&&
\bullet\lineto[r]&\overtag{\bullet}{-d}{8pt}\undertag{}{}{6pt}\lineto[r]&\bullet
\dashto[rr]&&\righttag{\bullet}{n_{t}/q_{t}}{8pt}\\&~\\&~\\&~\\&~}
$$

Here,  the continued fraction expansion $n/q=b_1-1/b_2-\cdots -1/b_s$ represents a string of rational curves emanating from the center:
\begin{equation}
\label{step1}
\xymatrix@R=6pt@C=24pt@M=0pt@W=0pt@H=0pt{
\\&\undertag{\bt}{-b_1}{6pt}\lineto[r]&\undertag{\bt}{-b_2}{6pt}
\dashto[r]&\dashto[r]
&\undertag{\bt}{-b_s}{6pt}\\
&&&&\\
&&&&\\
&&&&\\
}
\end{equation}
(We shall assume that $t\geq 3$.)   It is well-known that 
\begin{equation}
\label{step2}
\det \Gamma=\pm n_1n_2\cdots n_t(d-\sum_{i=1}^t(q_i/n_i)). 
\end{equation}
As long as $\det\Gamma\neq 0$, one can solve the equations $$K\cdot E_i + E_i\cdot E_i=-2,$$
and write $$K= \Sigma k_iE_i, \ \ \  k_i\in \mathbb Q.$$

Recall a negative-definite $\Gamma$ arises from the resolution of a weighted-homogeneous surface singularity.  In Section $2$ of \cite{w6}, the $k_i$ are computed in this case; but only non-degeneracy of $\Gamma$ was used, so the same formulas apply.

The formulas are expressed in terms of the two invariants
$$e=d-\sum_{i=1}^t(q_i/n_i) $$
$$\chi=t-2-\sum_{i=1}^t(1/n_i)=-2+\sum_{i=1}^t(1-1/n_i).$$
Since $e\neq 0$ by (2), we can define $\beta=\chi/e$.

As in \cite{w6}, we consider the rational cycle $-(K+E).$   For a cyclic quotient chain as in $(1)$, let  $F_1, \cdots ,F_s$ denote the curves.  Define the rational cycle $e_i$ by the property $e_i(F_j)=-\delta_{ij}$; it is effective (i.e., has strictly positive coefficients).  Then consider the cycle $Y=\beta e_1-e_s$  (even if $s=1$.) Denote by $Y_k$ the corresponding cycle for the $k$th string corresponding to $n_k/q_k$ in the graph of $\Gamma$, where $1\leq k \leq t.$  Denoting by $E_0$ the central curve of $\Gamma$, Proposition $2.3$ of \cite{w6} yields 
\begin{equation} \label{*}
-(K+E)=\sum_{k=1}^t  Y_k +\beta E_0.
\end{equation} 

\begin{lemma} Assume $\Gamma$ is one of the graphs in Table $A.2$.
\begin{enumerate}
\item $\chi \geq 0$, and $\chi=0$ exactly for the log-canonical singularities $\mathcal W(0,0,0)$, $\mathcal N(0,0,0)$, $\mathcal M(0,0,0).$
\item $e<0$.
\item $\beta<0$ in all cases except the three log-canonicals above, in which case it is $0$.
\item $|\beta|<1$.
\end{enumerate}
\begin{proof}  The first statement is a simple check (a sum of three reciprocals of integers is rarely at least $1$).  For the second, one need only consider the cases when $d=1$.  But all of those examples have $2$ chains of $-2$ curves emanating from the central curve; such a chain has $q=n-1$, so $q/n\geq 1/2$.  The third means checking that $\chi<|e|$, or $$t-2-\sum_{i=1}^t(1/n_i)<\sum_{i=1}^t(q_i/n_i)-d.$$  This statement turns out to be equivalent to that in \cite{w6}, Lemma 2.4; but in any case, it is an exercise.  (We quickly note that $d=-1$ in case $\mathcal W$; for $\mathcal N$, we have $d=0$ and some $q_i=n_i-1$; for type $\mathcal M$, two strings have $q_i=n_i-1$.) 
\end{proof}
\end{lemma}
\begin{proposition}  Suppose $(Z,E)$ is a surface of type $\Gamma$, where $\Gamma$ is a graph in Table $A.2$.  
\begin{enumerate}
\item The canonical divisor of $Z$ is $$K=\sum k_iE_i,$$
where for all $i$ $-1\leq k_i < 0$.
\item $k_i=-1$ only for the log-canonical cases  $\mathcal W(0,0,0)$, $\mathcal N(0,0,0)$, $\mathcal M(0,0,0),$ and then only at the central curve.
\end{enumerate}
\begin{proof}  Since the divisors $E_i$ span Pic $Z$ rationally, one can write the canonical divisor as $\sum k_iE_i.$  These coefficients can therefore be computed as above just from the graph.  In terms of the divisor $-(K+E)$, the claim is that its coefficients $-k_i-1$ are between $-1$ and $0$, and equal $0$ only at the center for the $3$ special cases.  Lemma 1.1 (3) verifies this assertion for the central curve.  

 It remains to show that the coefficients of $Y_k$ are strictly between $-1$ and $0$.  Writing $Y=\beta e_1-e_s$, all $e_i$ have strictly positive coefficients; as $\beta \leq 0$, the coefficients of $Y$ are strictly negative.  
 
 Next, writing $F=\sum F_j$, we claim that $$(F+Y)\cdot F_j \leq 0,\text{all}\ j;$$ this implies $F+Y$ has strictly positive coefficients, so the coefficients of $Y$ are bigger than $-1$ (one could not have $F=-Y$).  For $j=1$, the term in question is $1-b_1-\beta<2-b_1\leq 0$ (it does not matter if $s=1$).  An easier argument handles the other cases.
\end{proof}
\end{proposition}

\begin{corollary}  Let $(Z,E)$ be a surface of type $\Gamma$ as above, and $C\subset Z$ an irreducible curve with $C\cdot C <0$ and not a component of $E$.  Then 
\begin{enumerate}
\item $C$ is a smooth rational curve with $C\cdot C=-1$
\item If $C\cdot E_i=1$ for some $i$, then there is another $E_j$ with $C\cdot E_j>0$.
\end{enumerate}
\begin{proof}  $C\cdot E_j \geq 0$ for all $j$, and is positive for at least one $j$ because $E$ supports an ample divisor.  By Proposition 1.2, we have $K\cdot C<0$, so the usual adjunction formula yields that $C$ is a smooth rational $-1$ curve, and $C\cdot K=-1$.  In particular, $\sum (-k_j)C\cdot E_j=1$.  So, the second statement will follow once we exclude that for one of the three log-canonicals, there is a $-1$ curve which intersects the central curve transversally but does not intersect any other curve.  But in each of those cases, adding such a $-1$ curve to $E$ would give a non-degenerate curve configuration, so that its class could not be a rational combination of the components of $E$.
\end{proof}
\end{corollary}

We can paraphrase the last result by saying there are no ``free'' $-1$ curves, intersecting $E$ only once.

\section{How to find  sets of $-1$ curves}
If $(Z,E)$ is a $\Gamma$ surface, where $\Gamma$ is of type $\mathcal W, \mathcal N$, or $\mathcal M$, we wish to show that it is a Basic Model.  That means, $Z$ contains a set of three $-1$ curves which allow one to blow down to $\mathbb P^2$; the basic cases identify possible locations of the curves relative to $E$.   The blow downs give a projectively rigid configuration in $\mathbb P^2$ ($4$ lines in general position), from which the uniqueness of the $\Gamma$ surface follows.  

The method (initially analogous to the one used in \cite{bs}) is to blow up and down around the central curve of $E$ to produce $(Z',E')$, on which the new central curve $E'_0$ is a smooth rational curve of self-intersection $+1$.   The complete linear system associated to such a curve gives a birational map $\Phi:Z'\rightarrow \mathbb P^2$ which is an isomorphism in a neighborhood of $E'_0$.   It is the analysis of this map which will produce $-1$ curves first on $Z'$ and then on $Z$.  In each case, it will follow from the construction that one has an isomorphism of $Z'-E_0'$ with some open set in $Z$; we can conclude  (as in Corollary 1.3) that there are no ``free'' $-1$ curves on $Z'$, and a curve which is not a component of $E'$ has self-intersection $\geq -1$.

Here is how we proceed:

We note $\Phi(E'_0)\equiv L$ is a line. Each curve  $C$ in $E'$ adjacent to $E'_0$ is smooth (there are usually $3$ such), and $\Phi(C)$ is a (possibly singular) rational plane curve of degree $d=C\cdot E'_0 >0$.  The behavior of these image curves near $L$ is the same as it was on $Z'$, and the key will be to figure out their intersections away from $L$.  $\Phi(E')$ will have at most three singular points (away from $L$), and all possible configurations need to be considered.

The construction of $Z'$ shows that the components of $E'$ span Pic($Z'$) rationally, so no curves are disjoint from $E'$, and $\Phi$ is a sequence of blowing-up points over $\Phi(E')$ away from $L$. 
 
 We note $\Phi^{-1}(\Phi(E'))$ consists of $E'$ and (usually) three $-1$ curves, which is the same as for the basic cases (i.e., Basic Models).   For, the number of blow-downs given by $\Phi$ depends only on $K_{Z'}^2$, which is computed from $\Gamma'$, so is the same as in the basic case.  New curves added to $\Phi^{-1}(\Phi(E'))$ have negative self-intersection, so are $-1$ curves.

 Emanating from each adjacent curve $C$ is a chain (possibly empty) $\mathcal C(C)$ of rational curves, frequently with a long tail of $-2$ curves; it is disjoint from $E'_0$, so $\Phi$ sends it to a point $\Phi(\mathcal C(C))\in \Phi(C)$.  Being a smooth point of $\mathbb P^2$, its inverse image is a ``blow-downable configuration".  It contains $\mathcal C(C)$, at least one $-1$ curve, and any other chains $\mathcal C(\tilde{C})$ with the same image under $\Phi$.   The inverse image of a singular point of $\Phi(E')$ (of course, away from the line $L$) is either a union of chains and $-1$'s, or a single $-1$ curve intersecting only adjacent curves.  Since there are at most three new $-1$'s added in the inverse image of $\Phi(E')$, there are at most three singular points.
 
\begin{proposition}  $\Phi(\mathcal C(C))$ is a singular point of $\Phi(E')$.  More generally,  to go from $\mathbb P^2$ to $Z'$, one blows up only singular points of the inverse images of $\Phi(E')$.
\begin{proof}  Suppose a smooth point on a curve in the blow-up process is blown-up further.  Then the inverse image in $\Phi^{-1}(\Phi(E'))$ contains a smooth curve $C$ plus a blow-downable configuration attached transversally at a point of $C$.  This configuration has a $-1$ curve.  If it were at an end, this would be a free curve, a contradiction.  If not, it would be an interior curve, and removing it would leave a bunch of curves disjoint from $E'$. This also is a contradiction.\end{proof}
\end{proposition}

To unravel $\Phi$, one first examines the possible intersections of the images of the adjacent curves, noting that there are at most three singular points.  In each case,  $\Phi$ must factor via the minimal resolution of the singular points of $\Phi(E')$; one gathers  information about $\Phi^{-1}(\Phi(E'))$, such as possible valency of curves, or whether the $-1$ curves must intersect $E'$ transversally.  

Possible blow-downable configurations on $Z'$ are formed by putting together chains and $-1$'s.
There are limits to the location of $-1$ curves. 

\begin{remark} The following two configurations are not negative-definite:
$$
\xymatrix@R=6pt@C=24pt@M=0pt@W=0pt@H=0pt{
\\&&&&\\\
&&&&&\overtag{\bullet}{-1}{8pt}&&&&\\
&&&&&\lineto[u]&&&&\\
&&&&&\lineto[u]&&&&\\
&&&&&\lineto[u]&&&&\\
\undertag{\bullet}{-2}{4pt}\lineto[r]
&\undertag{\bullet}{-1}{4pt}\lineto[r]&\undertag{\bullet}{-2}{4pt}
&&\undertag{\bullet}{-2}{4pt}\lineto[r]
&\undertag{\bullet}{-2}{4pt}\lineto[r]\lineto[u]&\undertag{\bullet}{-2}{4pt}
&&\\
&&&&\\
&&&&\\
&&&&
}
$$
 \end{remark}                   
 The first example implies that two different chains cannot be connected at $-2$'s.  The second implies that a $-1$ curve could intersect a chain of $-2$ curves only at one of its ends.   But a connection at the beginning of a $-2$ chain (next to the adjacent curve) has consequences.
 
\begin{remark}
 The configuration 
 $$
\xymatrix@R=6pt@C=24pt@M=0pt@W=0pt@H=0pt{
\\&&&&\\
&&\overtag{\bullet}{C}{8pt}&&\\
&&\lineto[u]&\\
&&\bullet
\lineto[u]_(.2){-1}{}\\
&&\lineto[u]\\
&\undertag{\bt}{P}{4pt}\lineto[r]&\undertag{\bt}{-2}{4pt}\lineto[u]\lineto[r]&\undertag{\bullet}{-2}{4pt}&\\
&&&&\\
&&&&\\
}
$$ 
will, when blown-down, produce two curves which do not intersect transversally.  
\end{remark}
Consequently, if a chain emanating from the $P$ adjacent curve begins with two $-2$ curves, a $-1$ curve intersecting the first of these cannot intersect another curve, unless the final curve in $\mathbb P^2$ has a non-transversal intersection.  Variants of this situation will arise as well.

It is not a priori clear that the full inverse image of $\Phi(E')$ has normal crossings; while a $-1$ intersects transversally in a blow-downable configuration, it could in principle attach non-transversally to an adjacent curve or curves.

Finally, we introduce a notational convenience.  In the various blown-up spaces between $Z'$ and  $\mathbb P^2$, we shall frequently refer to the image of an adjacent curve $C$ as $C'$, or $C'(s)$ (when the self-intersection at that stage is $s$).




\section{Type $\Gamma=\mathcal W(p,q,r)$}
 $$
\xymatrix@R=6pt@C=24pt@M=0pt@W=0pt@H=0pt{
\\&&&&&&&&&&&\\\
&&&&&&(p+1)&&\\
&&&&&&{\hbox to 0pt{\hss$\overbrace{\hbox to 80pt{}}$\hss}}&&&&&&&\\
&&&&&&&&&&&&\\
&&&&&&&&&&&&\\
&&&&&\undertag{\bullet}{}{4pt}\dashto[r]&\dashto[r]&\undertag{\bullet}{}{4pt}&\\
&&&&&\lineto[u]\lineto[d]&&&&&&&&\\
&&&&&\lineto[u]&&&&&\\
&&(r+1)&&&\righttag{\bullet}{-(q+2)}{4pt}\lineto[u]&&&(q+1)\\
&&{\hbox to 0pt{\hss$\overbrace{\hbox to 80pt{}}$\hss}}&&&\lineto[u]&&&{\hbox to 0pt{\hss$\overbrace{\hbox to 80pt{}}$\hss}}&\\
&\undertag{\bullet}{}{4pt}\dashto[r]&\dashto[r]&\undertag{\bullet}{}{4pt}\lineto[r]&\undertag{\bullet}{-(p+2)}{4pt}\lineto[r]
&\undertag{\bullet}{+1}{4pt}\lineto[r]\lineto[u]&\undertag{\bullet}{-(r+2)}{4pt}\lineto[r]
&\undertag{\bullet}{}{4pt}\dashto[r]&\dashto[r]&\undertag{\bullet}{}{4pt}
&&\\
&&&&\\
&&&&\\
&&&&
}
$$
 Suppose we are given a surface $Z$ of type $\Gamma=\mathcal W(p,q,r)$.  The central curve $E_0$ has self-intersection $+1$, so in the above discussion we can set $Z'=Z$.  The three curves adjacent to $E_0$ are $P,Q,R$, with self-intersections respectively  $-(p+2), -(q+2), -(r+2)$, and with chains consisting solely of $-2$ curves.  We will prove the existence of the three rational $-1$ curves which appear in the Basic Model. (For the case $p=q=r$, there is a second choice, by flipping $Q$ and $R$ and their chains). Each $-1$ curve will connect an adjacent curve with the end of a chain associated to a different adjacent curve.  

By earlier discussion, $\Phi(P), \Phi(Q)$, and $\Phi(R)$ are lines intersecting $\Phi(E_0)=L$ in distinct points.  Thus $\Phi(E)$ either contains three lines through one point, or consists of four lines in general position.  But having a triple point would mean there is only one singular point, so that all three chains would be connected by three $-1$'s.  Remark 2.2 shows this is impossible.  So $\Phi(E)$ has three ordinary double points, and hence $\Phi^{-1}(\Phi(E))$ has normal crossings and only curves of valency two (of course, not counting intersection with $E_0$).

Each chain connects to other curves in $E$ only with a $-1$ attached at one of its ends.  That $-1$ cannot connect with another chain (Remark 2.2), so intersects an adjacent curve. The three $-1$'s are distributed among the three chains.   By Remark 2.3, the $-1$ curve appended to a chain must intersect at the far end. 

If the $-1$ at the end of $\mathcal C(P)$ intersects $Q$, then $\Phi(\mathcal C(P))$  is the intersection of the lines $\Phi(P)$ and $\Phi(Q)$.  Therefore, the intersection point $\Phi(P)\cap \Phi(R)$ must be the image of the chain $\mathcal C(R)$, whose $-1$ curve attachment at the end must intersect $P$.  Consequently, $\Phi(Q)\cap \Phi(R)$ comes from the chain $\mathcal C(Q)$, with a $-1$ attached to $R$.  Blowing-down $E$ completely in this way, collapsing first the $3$ $-1$'s and then the adjacent $-2$'s (which have become $-1$'s), one sees that this can happen only if $p=q=r$. 

On the other hand, if the $-1$ curve emanating from the end of $\mathcal C(P)$ intersects $R$, then the same analysis shows the resulting placement of two $-1$'s as before always blows down exactly to $4$ lines in general position.  Thus, one has a unique location of the $-1$'s (seen in Table $A.2$), except in case $p=q=r$, in which case there is a second possibility.   These are exactly the Basic Models for $\mathcal W$.

\section{ Type $\Gamma=\mathcal N(p,q,r)$, $p>0$}
 $$
\xymatrix@R=6pt@C=24pt@M=0pt@W=0pt@H=0pt{
\\&&&&&&&&&&&\\\
&&&&&&&(p+1)&&\\
&&&&&&&{\hbox to 0pt{\hss$\overbrace{\hbox to 80pt{}}$\hss}}&&&&&&&\\
&&&&&&\undertag{\bullet}{}{4pt}\dashto[r]&\dashto[r]&\undertag{\bullet}{}{4pt}&\\
&&&&&&\lineto[u]\lineto[d]&&&&&&&&\\
&&&&&&&&&&&&&&\\
&&(r+1)&&&&\lineto[u]&&&(q+2)\\
&&{\hbox to 0pt{\hss$\overbrace{\hbox to 80pt{}}$\hss}}&&&&\lineto[u]&&&{\hbox to 0pt{\hss$\overbrace{\hbox to 80pt{}}$\hss}}&\\
&\undertag{\bullet}{}{4pt}\dashto[r]&\dashto[r]&\undertag{\bullet}{}{4pt}\lineto[r]&\undertag{\bullet}{-(p+2)}{4pt}\lineto[r]&\undertag{\bullet}{-(q+2)}{4pt}\lineto[r]
&\undertag{\bullet}{0}{4pt}\lineto[r]\lineto[u]&\undertag{\bullet}{-(r+2)}{4pt}\lineto[r]
&\undertag{\bullet}{}{4pt}\dashto[r]&\dashto[r]&\undertag{\bullet}{}{4pt}
&&\\
&&&&\\
&&&&\\
&&&&
}
$$

\bigskip
We consider initially the case $p>0$.  Proceed to a new $(Z',E')$ as follows: First, blow-up any point on the central curve not on one of the three adjacent curves.  This makes the central curve a $-1$ curve, with four curves emanating from it, and adds a new curve $F$.  Now blow-down the old central curve and the one above it in the graph above, yielding:
 $$
\xymatrix@R=6pt@C=24pt@M=0pt@W=0pt@H=0pt{
\\&&&&&&&&&&&\\\
&&&&&&&&(p-1)&&\\
&&&&&&P&&{\hbox to 0pt{\hss$\overbrace{\hbox to 80pt{}}$\hss}}&&&&&&&\\
&&&&&&\lefttag{\bullet}{-1}{8pt}\lineto[r]&\undertag{\bullet}{}{4pt}\dashto[r]&\dashto[r]&\undertag{\bullet}{}{4pt}&\\
&&(r+1)&&&&\lineto[u]&&&(q+2)\\
&&{\hbox to 0pt{\hss$\overbrace{\hbox to 80pt{}}$\hss}}&&&Q&\lineto[u]&R&&{\hbox to 0pt{\hss$\overbrace{\hbox to 80pt{}}$\hss}}&\\
&\undertag{\bullet}{}{4pt}\dashto[r]&\dashto[r]&\undertag{\bullet}{}{4pt}\lineto[r]&\undertag{\bullet}{-(p+2)}{4pt}\lineto[r]&\undertag{\bullet}{-q}{4pt}\lineto[r]
&\undertag{\X}{+1}{8pt}\lineto[r]\lineto[u]&\undertag{\bullet}{-r}{4pt}\lineto[r]
&\undertag{\bullet}{}{4pt}\dashto[r]&\dashto[r]&\undertag{\bullet}{}{4pt}
&&\\
&&&&\\
&&&&\\
&&&&
}
$$
The curve  $F$ has now become a $+1$ central curve $E_0'$, intersecting transversally a $-1$ curve $P$, from which a chain of $(p-1)$ $-2$ curves emerge.  The two other original adjacent curves are still adjacent, but their self-intersections are now $-q$ and $-r$; we call the new ones $Q$ and $R$.  But now $Q$ and $R$ are simply tangent to each other and to $E_0'$, and $P,Q,R$ all intersect it at the same point of the central curve.  (We use the symbol $\X$ as a reminder that the intersections are not  transversal.)  The usual comparison with $Z$ shows that $Z'$ has no free $-1$'s, and the only negative curves off $E'$ are $-1$'s.

Proceeding as above, one constructs $\Phi:Z'\rightarrow \mathbb P^2$.  Then $\Phi(E')$ consists of the line $\Phi(E_0')=L$, two smooth conics $\Phi(Q)$ and $\Phi(R)$ intersecting each other and $L$ simply tangentially at a point of $L$, and a line $\Phi(P)$ intersecting transversally at that point.  By the usual argument, there are three additional $-1$'s in $\Phi^{-1}(\Phi(E'))$.   Here are the possibilities for the other intersections of the images of the three adjacent curves:
\begin{description}
\item [Case I] The two conics intersect tangentially at one other point, and the line passes through it (one singular point).
\item [Case II] The two conics intersect tangentially at one other point, and the line intersects each conic at a different point (three singular points).
\item [Case III] The two conics intersect transversally at two other points, and the line passes through one of these points (two singular points).
\end{description}
 
Resolving singularities in each case,  one finds that the inverse image of $\Phi(E')$ has normal crossings and a unique curve of valency three (of course, away from $E_0'$), which is not an adjacent curve.

 The inverse image of a singular point of $\Phi(E')$ is a blow-downable graph which is a combination of $-1$ curves and some of the three chains.  In particular, $\mathcal C(P)$ and $\mathcal C(R)$ each become blow-downable with a $-1$ curve appended at the end; further, each one could attach to a chain only at the $-(p+2)$ location of $\mathcal C(Q)$ (via Remark 2.2).  When $p=1$, then $\mathcal C(P)$ is the empty chain; but at least one $-1$ curve must still emerge from $P$, since $\Phi(P)$ intersects the other curves.
\subsection{Case I for $\mathcal N(p,q,r)$, $p>0$}
Case I does not occur.  Since one cannot have a valency four curve, a simple check shows there is no way to attach all three chains using three $-1$'s to get one blow-downable configuration (even in case $p=1$ when $\mathcal C(P)$ is empty).
\subsection{Case II for $\mathcal N(p,q,r)$, $p>0$}
In Case II, the conics $\Phi(Q)$ and $\Phi(R)$ are tangent away from $L$ and the line $\Phi(P)$ intersects each once, so there are three singular points.  That means $\mathcal C(Q)$ must become blow-downable either with the addition of a single $-1$ curve, or with a single $-1$ joining it and another chain.   One computes that adding a single $-1$ to $\mathcal C(Q)$ can make it blow-downable only if $p=r$ and the $-1$ is attached to the $-2$ curve adjacent to the $-(p+2)$.  That $-1$ curve must intersect one of the other adjacent curves.  But as in Remark 2.3, blowing down would give the image of that adjacent curve a worse than simple tangency with $\Phi(Q)$.  This is a contradiction. 

The only other option is to attach $\mathcal C(R)$ with a $-1$ adjoined at the $-(p+2)$ entry of $\mathcal C(Q)$; this blows down exactly when $p=q+2$.  However, if $r>0$ one sees that $\Phi(Q)$ and $\Phi(R)$ will have a higher order of tangency; this is ruled out.  

So, consider the special situation $p=q+2$ and $r=0$.  We show there is a unique way to find three $-1$'s which blow-down this $E'$.  The graph is
$$
\xymatrix@R=6pt@C=24pt@M=0pt@W=0pt@H=0pt{
\\&&&&&&&&&&&\\\
&&&&&&&&(q+1)&&\\
&&&&&&P&&{\hbox to 0pt{\hss$\overbrace{\hbox to 80pt{}}$\hss}}&&&&&&&\\
&&&&&&\lefttag{\bullet}{-1}{8pt}\lineto[r]&\undertag{\bullet}{}{4pt}\dashto[r]&\dashto[r]&\undertag{\bullet}{}{4pt}&\\
&&&&&&\lineto[u]&&&(q+2)\\
&&&&&Q&\lineto[u]&R&&{\hbox to 0pt{\hss$\overbrace{\hbox to 80pt{}}$\hss}}&\\
&&&\undertag{\bullet}{}{4pt}\lineto[r]&\undertag{\bullet}{-(q+4)}{4pt}\lineto[r]&\undertag{\bullet}{-q}{4pt}\lineto[r]
&\undertag{\X}{+1}{8pt}\lineto[r]\lineto[u]&\undertag{\bullet}{0}{4pt}\lineto[r]
&\undertag{\bullet}{}{4pt}\dashto[r]&\dashto[r]&\undertag{\bullet}{}{4pt}
&&\\
&&&&\\
&&&&\\
&&&&
}
$$
The above discussion states that the $-(q+4)$ curve has valency three in $\Phi^{-1}(\Phi(E'))$ and connects with a $-1$ curve from one of the ends of $\mathcal C(R)$.

Here is the minimal blow-up of $\Phi(E')$ which separates the line and two conics:
$$
\xymatrix@R=6pt@C=24pt@M=0pt@W=0pt@H=0pt{
\\
&&&&&&&&&&&\\
&&&\righttag{\bullet}{M'(-2)}{4pt}\lineto[d]&&&&\\
&&&\lineto[d]&&&&\\
&&&\lineto[d]&&&&\\
&&\lefttag{\bullet}{R'(1)}{4pt}\lineto[r]&\undertag{\bullet}{N'(-1)}{2pt}\lineto[u]\lineto[r]&\righttag{\bullet}{Q'(1)}{4pt}&&&&&&\\
&&\lineto[u]&&\lineto[u]&&\\
&&\lineto[u]&&\lineto[u]&&\\
&&\lineto[u]&&\lineto[u]&&\\
&&\undertag{\bullet}{-1}{4pt}\lineto[r]\lineto[u]&\undertag{\bullet}{P'(-1)}{4pt}\lineto[r]&\undertag{\bullet}{-1}{4pt}\lineto[u]&&
&&\\
&&&&\\
&&&&\\
&&&&
}
$$
Recall  $P'(-1), Q'(1), R'(1)$ are the images of the adjacent curves on partial blow-ups plus their self-intersections there.  New curves $M'$ and $N'$ have been named. We specify what is needed in order to blow-up further to get to $E'$.
Since $P$ has degree $-1$, there can be no further blow-ups along the bottom line.  One cannot blow-up between $N'(-1)$ and $M'(-2)$, as an $M'(-3)$ (with degree $\leq -3$) would eventually become the $-(q+4)$ curve, but not adjacent to $Q$.  $N'(-1)$ must be blown up somewhere, else the inverse image of a singular point would be a single $-1$ and $-2$, but not intersecting $P$.   Therefore, $N'$ is the valency three curve which will eventually become the $-(q+4)$ curve.  But that curve is adjacent to $Q$, so the only place $N'(-1)$ can be blown-up is at the intersection with $R'(1)$. After one blow-up, one reaches 
$$
\xymatrix@R=6pt@C=24pt@M=0pt@W=0pt@H=0pt{
\\&&&&\\\
&&&&\\
&\dashto[r]&\undertag{\bullet}{R'(0)}{4pt}\lineto[r]
&\undertag{\bullet}{-1}{4pt}\lineto[r]&\undertag{\bullet}{N'(-2)}{4pt}\dashto[r]\lineto[u]
&&\\
&&&&\\
&&&&\\
&&&&
}
$$
As $R$ has self-intersection $0$, the only further blowing-up takes place between the $-1$'s and the curve on the right.  After $q+2$ of such blow-ups, one has $R'(0)$ followed by $(q+2)$ $-2$'s, followed by a $-1$, followed by $N'(-(q+4))$. This completes most of the blow-up to $E'$.   It remains only to complete the portion
$$
\xymatrix@R=6pt@C=24pt@M=0pt@W=0pt@H=0pt{
\\&&&&\\\
&\dashto[r]&\undertag{\bullet}{P'(-1))}{4pt}\lineto[r]
&\undertag{\bullet}{-1}{4pt}\lineto[r]&\undertag{\bullet}{Q'(1)}{4pt}\dashto[r]
&&\\
&&&&\\
&&&&\\
&&&&
}
$$
This must be done in the usual way of blowing up $(q+1)$ times between the $-1$ and $Q'(1)$.  This completes the blowing up to reach $E'$.  

We note the locations of the three $-1$'s which allow the blow-down: between the end of $\mathcal C(R)$ and the $-(q+4)$ curve; between the end of $\mathcal C(P)$ and $Q$; and between $P$ and $R$.  As mentioned before, pulling these curves back to $Z$ gives two of the $-1$'s seen by using the symmetry (given the special values of $p,q,r$) between the top and right hand chains in the graph.  The third of the $-1$'s on $Z$ can be found in the above example by pulling back from $Z'$ the inverse image of one of the lines through an intersection of the two conics.

\subsection{Case III for $\mathcal N(p,q,r)$, $p>0$}  There are two singular points; the two conics meet in two distinct points, and the line $\Phi(P)$ passes through one of those.  The minimal resolution is
$$
\xymatrix@R=6pt@C=24pt@M=0pt@W=0pt@H=0pt{
\\&&&&&&&&&\\\
&&&&&&\overtag{\bullet}{Q'(2)}{12pt}\lineto[r]&\lineto[r]&\overtag{\bullet}{-1}{12pt}&\\
&&&&&&\lineto[u]&&\lineto[u]&&&&&\\
&&&&&&\lineto[u]&&\lineto[u]&&&&&\\
&&&&\undertag{\bullet}{P'(0)}{4pt}\lineto[r]&\lineto[r]&\undertag{\bullet}{M'(-1)}{4pt}\lineto[r]\lineto[u]&\lineto[r]&\undertag{\bullet}{R'(2)}{4pt}\lineto[u]&&\\
&&&&\\
&&&&\\
&&&&
}
$$
As $P'(0)$ has but one connection with the rest of the graph, as before one must blow up between it and $M'(-1)$, eventually reaching $P'(-1)$ followed by $(p-1)$ $-2$'s, followed by a $-1$, followed by $M'(-(p+1))$.  If there were no further blow-ups between $M'(-(p+1))$ and $Q'(2)$ or $R'(2)$, then $E'$ itself would contain a curve intersecting both $R$ and $Q$.  This does not happen, so the $M'$ curve must be blown-up at least once more and become $M'(-(p+2))$, the $-(p+2)$ curve in $\mathcal C(Q)$ adjacent to $Q$.

Thus, one has $P$ followed by $\mathcal C(P)$ (even if empty) followed by $-1$ attached to the $-(p+2)$ curve.   Since $Q$ is adjacent to that curve, in the above diagram there is no blowing-up between them; to reach $Q$ of degree $-q$, one blows up repeatedly along the top row.  This yields $Q$ followed by a $-1$ followed by $(q+2)$ $-2$'s followed by $R'(2)$.  In other words, $Q$ is attached via a $-1$ with the end of $\mathcal C(R)$, accounting for the intersection point of the two conics not involving the line.   The only way for $\mathcal C(Q)$ to attach is via a $-1$ at its end intersecting $R$.

This is exactly the placement of the three $-1$'s on $Z'$ as would happen in the basic case.  Passing from $Z'$ back to $Z$, the $-1$'s and their relative location stays the same.  We recover a Basic Model.

\subsection{Type $\Gamma=\mathcal N(0,q,r)$}

In this situation, we proceed from $(Z,E)$ to $(Z',E')$ exactly as before, first by blowing up to add one curve, and then blowing down two curves.  The difference is that there is now no longer an adjacent curve $P$; rather, there are just two adjacent curves $Q$ and $R$, each simply tangent to each other and to the central curve $E'_0$.  In this case, the map $\Phi$ sends $E'$ to the line $L$ and the two conics $\Phi(Q)$ and $\Phi(R)$, and two $-1$'s are needed to make the blow-down.

As $p=0$, the chains $\mathcal C(Q)$ and $\mathcal C(R)$ consist solely of $-2$-curves, so cannot be connected via a $-1$.  They are individually blow-downable by the usual addition of a $-1$ curve at the beginning or end of the chain.  So, $\Phi(E')$ has two singular points, hence the conics intersect transversally.  If either of these $-1$'s occurred at the beginning of a chain, then blowing down the chain would give a tangency or worse.   Thus, the $-1$'s are at the far ends of the chains, and blowing down each chain gives one of the two intersection points of the conics $\Phi(Q)$ and $\Phi(R)$.  Thus, there is a unique blow-down, so the basic case is the only one.

If one pulls these two $-1$'s from $Z'$ back to $Z$, one might ask where is the third $-1$ needed for blow-down.   This can be found by pulling back to $Z'$ and then $Z$ one of the lines connecting the central point of $L$ with an intersection point of the two quadrics.  Again, this is a Basic Model.

\section{ Type $\Gamma=\mathcal M(p,q,r)$}
$$
\xymatrix@R=6pt@C=24pt@M=0pt@W=0pt@H=0pt{
\\&&&&&&&&&&&\\\
&&&&&(p+1)&&\\
&&&&&{\hbox to 0pt{\hss$\overbrace{\hbox to 80pt{}}$\hss}}&&&&&&&\\
&&&&\undertag{\bullet}{}{4pt}\dashto[r]&\dashto[r]&\undertag{\bullet}{}{4pt}&\\
&&&&\lineto[u]\lineto[d]&&&&&&&&\\
&&&&&&&&&&&&\\
&&(r+2)&&\lineto[u]&&&&&(q+2)\\
&&{\hbox to 0pt{\hss$\overbrace{\hbox to 80pt{}}$\hss}}&&\lineto[u]&&&&&{\hbox to 0pt{\hss$\overbrace{\hbox to 80pt{}}$\hss}}&\\
&\undertag{\bullet}{}{4pt}\dashto[r]&\dashto[r]&\undertag{\bullet}{}{4pt}\lineto[r]&\undertag{\bullet}{-1}{4pt}\lineto[r]\lineto[u]&\undertag{\bullet}{-(q+2)}{4pt}\lineto[r]
&\undertag{\bullet}{-(r+2)}{4pt}\lineto[r]&\undertag{\bullet}{-(p+2)}{4pt}\lineto[r]
&\undertag{\bullet}{}{4pt}\dashto[r]&\dashto[r]&\undertag{\bullet}{}{4pt}
&&\\
&&&&\\
&&&&\\
&&&&
}
$$

We refer to the three strings emanating from the central $-1$ curve in the diagram as the $r, p, $ and $q$ directions.  Start with the assumption that $p,r>0.$  To form the desired $(Z',E')$,  blow-down the central $-1$ curve and the first two curves in the $r$ direction.  Then the first curve in the $p$ direction is the new $E'_0$, of self-intersection $+1$, and $E'$ is 
$$
\xymatrix@R=6pt@C=24pt@M=0pt@W=0pt@H=0pt{
\\&&&&&&&&&&&\\\ 
&&&&&&&(p-1)&&\\
&&&&&&&{\hbox to 0pt{\hss$\overbrace{\hbox to 80pt{}}$\hss}}&&&&&&&\\
&&&&&\overtag{\bullet}{P}{12pt}\lineto[r]&\undertag{\bullet}{}{4pt}\dashto[r]&\dashto[r]&\undertag{\bullet}{}{4pt}&\\
&&(r-1)&&&\lineto[u]&&&&&(q+2)\\
&&{\hbox to 0pt{\hss$\overbrace{\hbox to 80pt{}}$\hss}}&&R&\lineto[u]&Q&&&&{\hbox to 0pt{\hss$\overbrace{\hbox to 80pt{}}$\hss}}&\\
&\undertag{\bullet}{}{4pt}\dashto[r]&\dashto[r]&\undertag{\bullet}{}{4pt}\lineto[r]&\undertag{\bullet}{-1}{4pt}\lineto[r]&\undertag{\X}{+1}{8pt}\lineto[r]\lineto[u]
&\undertag{\bullet}{-(q-1)}{4pt}\lineto[r]&\undertag{\bullet}{-(r+2)}{4pt}\lineto[r]
&\undertag{\bullet}{-(p+2)}{4pt}\lineto[r]
&\undertag{\bullet}{}{4pt}\dashto[r]&\dashto[r]&\undertag{\bullet}{}{4pt}
&&\\
&&&&\\
&&&&\\
&&&&
}
$$

 The new central configuration consists of a $+1$ curve $E'_0$ and three adjacent curves $P,Q,R$, but their intersections are no longer transversal.  $Q$ is a $-(q-1)$ curve with a tangency of order $3$ with $E'_0$ at a point, and $R$ intersects $E'_0$ at that point, transversally to both $E'_0$ and $Q$.  Finally, $P$ intersects $E'_0$ transversally at a different point.  Note that again there are no free $-1$ curves on $Z'$, given this property on $Z$.

The new map $\Phi:Z'\rightarrow \mathbb P^2$ arising from $E'_0$ makes $\Phi(E'_0)=L$ a line, $\Phi(Q)$ a rational cubic curve with triple tangency at a point of $L$, $\Phi(R)$ a line through that central point transversal to both $L$ and the cubic, and $\Phi(P)$ a line transversal to $L$ intersecting at a different point.  

A rational cubic curve is either nodal or cuspidal, and each type is unique up to projective equivalence.  Each has a unique flex point, i.e., smooth point whose tangent line intersects with multiplicity $3$. A calculation shows that a line through that flex point cannot be tangent to the curve at a smooth point.  As a result, here are the only possible configurations of $\Phi(E')$ with at most three singular points:
\begin{description}
\item [Case $I_c$ (resp. $I_n$)] Cubic is cuspidal (resp. nodal), $\Phi(R)$ passes through the singular point, $\Phi(P)$ has multiplicity three at that point (one singular point for $\Phi(E')$)
\item [Case $II_c$ (resp. $II_n$)] Cubic is cuspidal (resp. nodal), $\Phi(R)$ passes through the singular point, $\Phi(P)$ passes through the singular point plus another point of the cubic (two singular points)
\item [Case $III_c$ (resp. $III_n$)] Cubic is cuspidal (resp. nodal), $\Phi(R)$ intersects the cubic in two distinct smooth points, $\Phi(P)$ passes through one of those two points and the singular point (three singular points)
\end{description}

We rule out all but Case $I_n$  using the interplay between requirements of the graph $E'$ and resolution of the singularities of $\Phi(E')$. One cannot dismiss a priori the occurrence of a non-transversal intersection of a $-1$ curve with another curve.  

\subsection{Case $III$ for $\mathcal M(p,q,r)$, $p,r>0$} In Case $III$, $\Phi(E')$ would have an ordinary triple point involving all three curves.  The inverse image of this singular point would be a single $-1$ curve with valency (at least) three, intersecting all adjacent curves or their chains.  We show this cannot happen.  Blowing up the three singular points yields the curve configuration

$$
\xymatrix@R=6pt@C=24pt@M=0pt@W=0pt@H=0pt{
\\&&&&&&&&&\\\
&&&&\overtag{\bullet}{-1}{12pt}\dashto[r]&\dashto[r]&\overtag{\bullet}{Q'}{12pt}\lineto[r]&\lineto[r]&\overtag{\bullet}{-1}{12pt}&\\
&&&&\lineto[u]&&\lineto[u]&&\lineto[u]&&&&&\\
&&&&\lineto[u]&&\lineto[u]&&\lineto[u]&&&&&\\
&&&&\undertag{\bullet}{P'(-1)}{4pt}\lineto[u]\lineto[r]&\lineto[r]&\undertag{\bullet}{-1}{4pt}\lineto[r]\lineto[u]&\lineto[r]&\undertag{\bullet}{R'(-1)}{4pt}\lineto[u]&&\\
&&&&\\
&&&&\\
&&&&
}
$$
Here as before $P'(-1)$ indicates the image of $P$, with self-intersection $-1$.   The dotted connection on the top row between a $-1$ and $Q'$ indicates that the intersection of the curves is not transversal; $Q'$ intersects tangentially in the cusp case, and in two points in the nodal.   Since $P$ has self-intersection $-2$, to reach $Z'$ one needs to blow-up the intersection of $P'(-1)$ with exactly one of its neighboring $-1$'s.  If one blows up at the $-1$ above, that would make $Q'$ intersect non-transversally with a $-2$ curve, which must be resolved.  That further resolution would introduce a fourth $-1$ curve disjoint from the other three; this is a contradiction.  

So, one would have to blow up between $P'(-1)$ and the $-1$ on its right, converting that valency three curve to a $-2$, hence no longer eligible to be the $-1$ of valency three or more.  The only other possible way to get a trivalent $-1$ would be to blow up the non-transversal intersection between $Q'$ and the $-1$ on its left.  In the case of a node, this would result in a new trivalent curve, but it would be a $-2$.  In the case of a cusp, the only way to get a trivalent curve would be to blow-up twice, in which case the original $-1$ intersecting $Q'$ would become a $-3$.  However, $P$ does not intersect a curve of degree $\leq -3$ (only a $-2$ or $-1$ is allowed).  This completes the argument.

\subsection{Case $II$ for $\mathcal M(p,q,r)$, $p,r>0$} $\Phi(E')$ has two intersection points, with minimal resolution 
$$
\xymatrix@R=6pt@C=24pt@M=0pt@W=0pt@H=0pt{
\\&&&&&&&&&\\\
&&&&\overtag{\bullet}{Q'}{12pt}\dashto[r]&\dashto[r]&\overtag{\bullet}{-1}{12pt}\lineto[r]&\lineto[r]&\overtag{\bullet}{R'(0)}{12pt}&\\
&&&&\lineto[u]&&\lineto[u]&&&&&&&\\
&&&&\lineto[u]&&\lineto[u]&&&&&&&\\
&&&&\undertag{\bullet}{-1}{4pt}\lineto[u]\lineto[r]&\lineto[r]&\undertag{\bullet}{P'(0)}{4pt}\lineto[u]&&\\
&&&&\\
&&&&\\
&&&&
}
$$
The dotted line again refers to the non-transverse intersections due to the cusp or node.  Since $R$ has self-intersection $-1$, one must blow up the point between $R'(0)$ and its neighbor, converting it to a $-2$ and obtaining a new $-1$.  But one cannot have a non-transverse intersection between $Q'$ and a $-2$, so that must be blown up, converting the $-2$ to a $-3$ and inserting a $-1$ connecting it to $Q'$ (in both the node and cusp cases). That makes three disjoint $-1$'s.  But $P'(0)$ now intersects a $-3$ curve, yet $P$ does not; so, that intersection point must be blown up, giving a fourth $-1$, which is not allowed.

\subsection{Case $I$ for $\mathcal M(p,q,r)$, $p,r>0$}
Here $\Phi(E')$ has one singular point.  Blowing it up gives a $-1$ curve we shall call $M'(-1)$.  It has a simple intersection with $R'(0)$, which intersects no other curves. Since $R$ has self-intersection $-1$, the only way to achieve that is to blow-up the intersection, giving
$$
\xymatrix@R=6pt@C=24pt@M=0pt@W=0pt@H=0pt{
\\&&&&&&&&&&\\
&\undertag{\bullet}{R'(-1)}{4pt}\lineto[r]&\undertag{\bullet}{-1}{4pt}\lineto[r]&
\undertag{\bullet}{M'(-2)}{4pt}\dashto[r]&
&&\\
&&&&\\
&&&&\\
&&&&
}
$$
Further, the only way to achieve the string $\mathcal C(R)$ with $(r-1)$ $-2$ curves is to repeatedly blow up $-1$'s away from $R'(-1)$:
$$
\xymatrix@R=6pt@C=24pt@M=0pt@W=0pt@H=0pt{
\\&&&&&&&&&&&\\\
&&&&(r-1)&&&&&&&&\\
&&&&{\hbox to 0pt{\hss$\overbrace{\hbox to 80pt{}}$\hss}}&&&&&&\\
&&\undertag{\bullet}{R'(-1)}{4pt}\lineto[r]&\undertag{\bullet}{}{4pt}\dashto[r]&\dashto[r]&\undertag{\bullet}{}{4pt}\lineto[r]&\undertag{\bullet}{-1}{4pt}\lineto[r]&\lineto[r]&\undertag{\bullet}{M'(-(r+1))}{4pt}\dashto[r]
&&&&\\
&&&&\\
&&&&\\
&&&&
}
$$

Consider next the other intersections of $M'$. The cubic curve has become $Q'(5)$.  

In the cuspidal case, $M'(-(r+1))$ intersects $Q'(5)$ tangentially at a point through which $P'(0)$ passes transversally.  Blowing up that point gives a new $-1$ curve $N'(-1)$.  It has a simple intersection with $P'(-1)$, which intersects no other curves.  Since $P$ has self-intersection $-2$, as above one has to blow-up this intersection point.  Continuing as before, one gets part of the graph around $N'$ as
$$
\xymatrix@R=6pt@C=24pt@M=0pt@W=0pt@H=0pt{
\\&&&&&&&&&&&\\\
&&&&(p-1)&&&&&&&&\\
&&&&{\hbox to 0pt{\hss$\overbrace{\hbox to 80pt{}}$\hss}}&&&&&&\\
&&\undertag{\bullet}{P'(-2)}{4pt}\lineto[r]&\undertag{\bullet}{}{4pt}\dashto[r]&\dashto[r]&\undertag{\bullet}{}{4pt}\lineto[r]&\undertag{\bullet}{-1}{4pt}\lineto[r]&\lineto[r]&\undertag{\bullet}{N'(-(p+1))}{4pt}\dashto[r]
&&&&\\
&&&&\\
&&&&\\
&&&&
}
$$
There is as well a point where $N'$, $M'$, and $Q'$ intersect transversally.   Since all these curves will appear in $E'$, this point must be blown-up.  Putting everything together gives the graph
$$
\xymatrix@R=6pt@C=24pt@M=0pt@W=0pt@H=0pt{
\\&&&&&&&&&&&\\\
&&(r-1)&&&&\overtag{\bullet}{Q'}{8pt}&&&&(p-1)&\\
&&{\hbox to 0pt{\hss$\overbrace{\hbox to 80pt{}}$\hss}}&&&&\lineto[u]&&&&{\hbox to 0pt{\hss$\overbrace{\hbox to 80pt{}}$\hss}}&\\
\undertag{\bullet}{R'(-1)}{4pt}\lineto[r]&\undertag{\bullet}{}{4pt}\dashto[r]&\dashto[r]&\undertag{\bullet}{}{4pt}\lineto[r]&\overtag{\bullet}{-1}{8pt}\lineto[r]&\undertag{\bullet}{M'(-(r+3))}{4pt}\lineto[r]&\undertag{\bullet}{-1}{3pt}\lineto[r]\lineto[u]&\undertag{\bullet}{N'(-(p+2))}{4pt}\lineto[r]&\overtag{\bullet}{-1}{8pt}\lineto[r]&\undertag{\bullet}{}{4pt}\dashto[r]&\dashto[r]&\undertag{\bullet}{}{4pt}\lineto[r]&\undertag{\bullet}{P'(-2)}{4pt}
&&&&\\
&&&&\\
&&&&\\
&&&&
}
$$
But this configuration cannot be completed to produce $E'$ plus three $-1$ curves.  For instance, the $-1$ appended to the chain $\mathcal C(R)$ would intersect a curve whose valency remains two, hence an end of a chain.  But the self-intersection of that curve is $\leq -3$, an impossibility.  The cuspidal case is eliminated.

 We are down to the nodal case.  $M'(-(r+1))$ intersects $Q'(5)$ transversally in two points, through one of which $P'(0)$ passes with a third tangent direction.  Following the same procedure as in the cuspidal case, one reaches the graph
 $$
\xymatrix@R=6pt@C=24pt@M=0pt@W=0pt@H=0pt{
\\&&&&&&&&&&&\\\
&&(r-1)&&&&&&&&(p-1)&\\
&&{\hbox to 0pt{\hss$\overbrace{\hbox to 80pt{}}$\hss}}&&&&\overtag{\bullet}{Q'(4)}{10pt}&&&&{\hbox to 0pt{\hss$\overbrace{\hbox to 80pt{}}$\hss}}&\\
\undertag{\bullet}{R'(-1)}{4pt}\lineto[r]&\undertag{\bullet}{}{4pt}\dashto[r]&\dashto[r]&\undertag{\bullet}{}{4pt}\lineto[r]&\overtag{\bullet}{-1}{8pt}\lineto[r]&\undertag{\bullet}{M'(-(r+2))}{4pt}\lineto[r]\lineto[ur]&\lineto[r]&\undertag{\bullet}{N'(-(p+1))}{4pt}\lineto[r]\lineto[ul]&\overtag{\bullet}{-1}{8pt}\lineto[r]&\undertag{\bullet}{}{4pt}\dashto[r]&\dashto[r]&\undertag{\bullet}{}{4pt}\lineto[r]&\undertag{\bullet}{P'(-2)}{4pt}
&&&&\\
&&&&\\
&&&&\\
&&&&
}
$$
Now, $Q$ has self-intersection  $-(q-1)$, so $q+3$ further blow-ups are needed next to $Q'(4)$. 

Suppose the first blow-up is between it and $M'(-(r+2))$. Then the new $-1$ is the third one in the diagram, so all future blow-ups of $Q'$ must be adjacent to a $-1$.  The final position between $M'$ and $Q'$ is therefore
$$
\xymatrix@R=6pt@C=24pt@M=0pt@W=0pt@H=0pt{
\\&&&&&&&&&&&\\\
&&&&(q+2)&&&&&&&&\\
&&&&{\hbox to 0pt{\hss$\overbrace{\hbox to 80pt{}}$\hss}}&&&&&&\\
&\dashto[r]&\undertag{\bullet}{M'(-(r+3))}{4pt}\lineto[r]&\undertag{\bullet}{}{4pt}\dashto[r]&\dashto[r]&\undertag{\bullet}{}{4pt}\lineto[r]&\undertag{\bullet}{-1}{4pt}\lineto[r]&\lineto[r]&\undertag{\bullet}{Q'(-(q-1))}{4pt}\dashto[r]
&&&&\\
&&&&\\
&&&&\\
&&&&
}
$$
Note that $N'(-(p+1))$ still intersects $M'$ and $Q'$. This fits with the original $E'$ exactly when $r+2=p+1$ and $r+3=p+2$, i.e., $p=r+1$.  In that situation, the three desired $-1$'s connect the ends of the three chains as follows:  $\mathcal C(P)$ to the $-(r+2)$ curve; $\mathcal C(R)$ to the $-(p+2)$ curve; and $\mathcal C(Q)$ to $Q$.  Pulling these $-1$'s back to $Z$ gives the Basic Model for the special value $p=r+1$, when the graph has an obvious symmetry.

If the first blow-up takes place between $Q'(4)$ and $N'(-(p+1))$, then the same procedure as above gives $E'$ for all values of $p, q, r$.  Here, the $-1$ locations with the ends of chains are: $\mathcal C(P)$ to the $-(p+2)$ curve; $\mathcal C(R)$ to the $-(r+2)$ curve; and $\mathcal C(Q)$ to $Q$.  Pulling back to $Z$, one recovers the $-1$'s for the Basic Model, for all values $p, q, r$ (assuming still that $p,r>0)$.
  

\subsection{Type $\Gamma=\mathcal M(0,q,r), r \geq 1$}

We construct $Z'$ as above, noting that the one curve in the $p$ direction is the new central curve.  The diagram is as before, with only two curves adjacent to the central $E'_0$:
$$
\xymatrix@R=6pt@C=24pt@M=0pt@W=0pt@H=0pt{
\\&&&&&&&&&&&\\\
&&(r-1)&&&&&&&&(q+2)\\
&&{\hbox to 0pt{\hss$\overbrace{\hbox to 80pt{}}$\hss}}&&R&&Q&&&&{\hbox to 0pt{\hss$\overbrace{\hbox to 80pt{}}$\hss}}&\\
&\undertag{\bullet}{}{4pt}\dashto[r]&\dashto[r]&\undertag{\bullet}{}{4pt}\lineto[r]&\undertag{\bullet}{-1}{4pt}\lineto[r]&\undertag{\Circ}{+1}{8pt}\lineto[r]
&\undertag{\bullet}{-(q-1)}{4pt}\lineto[r]&\undertag{\bullet}{-(r+2)}{4pt}\lineto[r]
&\undertag{\bullet}{-2}{4pt}\lineto[r]
&\undertag{\bullet}{}{4pt}\dashto[r]&\dashto[r]&\undertag{\bullet}{}{4pt}
&&\\
&&&&\\
&&&&\\
&&&&
}
$$
As before, $Q$ is a $-(q-1)$ curve with a tangency of order three with $E_0'$ at a point, and $R$ intersects $E_0'$ at that point, transversally.  The basic case requires two $-1$'s to blow down, so the same should be true in general.

If  the line $\Phi(R)$ did not pass through the singular point of the cubic $\Phi(Q)$, it would intersect it in two distinct points, giving three singular points on $\Phi(E')$.  This cannot happen, so $\Phi(R)$ passes through the singular point.  If the singular point is a cusp, the same argument as above produces the same contradiction: the $-1$ at the end of $\mathcal C(R)$ would be intersecting a curve of valency two and self-intersection less than or equal to $-3$.  

So, the singular point is a node, and one gets a picture as before which gives a $-1$ at the end of $\mathcal C(R)$:
$$
\xymatrix@R=6pt@C=24pt@M=0pt@W=0pt@H=0pt{
\\&&&&&&&&&&&\\\
&&(r-1)&&&&&&&&&\\
&&{\hbox to 0pt{\hss$\overbrace{\hbox to 80pt{}}$\hss}}&&&&\overtag{\bullet}{Q'(4)}{10pt}&&&&\\
\undertag{\bullet}{R'(-1)}{4pt}\lineto[r]&\undertag{\bullet}{}{4pt}\dashto[r]&\dashto[r]&\undertag{\bullet}{}{4pt}\lineto[r]&\overtag{\bullet}{-1}{8pt}\lineto[r]&\undertag{\bullet}{M'(-(r+2))}{4pt}\lineto[r]\lineto[ur]&\lineto[r]&\undertag{\bullet}{-1}{4pt}\lineto[ul]&&&&&\\
&&&&\\
&&&&\\
&&&&
}
$$
As before, there is a unique way to complete to $E'$, with $-1$ curves appended to the end of chains as follows: $\mathcal C(R)$ at the $-(r+2)$ curve, and $\mathcal C(Q)$ at $Q$.  Pulling back to $Z$, and adding on the pull-back of the line on $Z'$ through the node and the central point of $L$, gives the Basic Model.


\subsection{Type $\Gamma=\mathcal M(p,q,0), p\geq 1$}

Again, we have the same $Z'$, but $E'$ is now
$$
\xymatrix@R=6pt@C=24pt@M=0pt@W=0pt@H=0pt{
\\&&&&&&&&&&&\\\
&&(p-1)&&&&&&&&(q+2)\\
&&{\hbox to 0pt{\hss$\overbrace{\hbox to 80pt{}}$\hss}}&&P&&Q&&&&{\hbox to 0pt{\hss$\overbrace{\hbox to 80pt{}}$\hss}}&\\
&\undertag{\bullet}{}{4pt}\dashto[r]&\dashto[r]&\undertag{\bullet}{}{4pt}\lineto[r]&\undertag{\bullet}{}{4pt}\lineto[r]&\undertag{\Circ}{+1}{8pt}\lineto[r]
&\undertag{\bullet}{-(q-1)}{4pt}\lineto[r]&\undertag{\bullet}{-2}{4pt}\lineto[r]
&\undertag{\bullet}{-(p+2)}{4pt}\lineto[r]
&\undertag{\bullet}{}{4pt}\dashto[r]&\dashto[r]&\undertag{\bullet}{}{4pt}
&&\\
&&&&\\
&&&&\\
&&&&
}
$$
Again $Q$ intersects $E_0'$ with multiplicity three, while $P$ intersects it transversally at a different point.  As in the basic case, two $-1$'s are needed to blow-down.

We rule out that $\Phi(E')$ could have two singular points, for which $\Phi(P)$ intersects the cubic at the singular point plus one other.  Note that $\mathcal C(Q)$ can be made blow-downable with the addition of a single $-1$ only in case $p=q+2$, and the $-1$ is added to the $-2$ to the right of the $-(p+2)$.  In that case, one would have a $-2$ curve of valency three.  A brief check of the resolution of $\Phi(E')$ shows that a valency three curve arises only when all non-transversality (from the cusp or node) has been resolved; thus, $\Phi^{-1}(\Phi(E'))$ has normal crossings.  Therefore the appended $-1$ curve to $\mathcal C(Q)$ intersects $P$ or $Q$ transversally.  But blowing-down produces a higher order tangency between two curves, as in Remark $2.3$.  So this case is eliminated.  In addition, one cannot combine $\mathcal C(Q)$ and $\mathcal C(P)$ with one $-1$ into a blow-downable configuration.

That leaves the case of one singular point of $\Phi(E')$, with $\Phi(P)$ intersecting there with multiplicity three.  The same argument as in the original situation shows one cannot have a cuspidal cubic, else $\mathcal C(P)$ (or just $P$, if $p=1$) would have a $-1$ intersecting a valency two curve with self-intersection $\leq -3$.  So, one must have a nodal cubic, for which $\Phi(P)$ is tangent to one of the branches of the node.  The same partial resolution as above is
$$
\xymatrix@R=6pt@C=24pt@M=0pt@W=0pt@H=0pt{
\\&&&&&&&&&&&\\\
&&&&&&&&&&(p-1)&\\
&&&&&&\overtag{\bullet}{Q'(4)}{10pt}&&&&{\hbox to 0pt{\hss$\overbrace{\hbox to 80pt{}}$\hss}}&\\
&&&&&\undertag{\bullet}{M'(-2)}{4pt}\lineto[r]\lineto[ur]&\lineto[r]&\undertag{\bullet}{N'(-(p+1))}{4pt}\lineto[r]\lineto[ul]&\overtag{\bullet}{-1}{8pt}\lineto[r]&\undertag{\bullet}{}{4pt}\dashto[r]&\dashto[r]&\undertag{\bullet}{}{4pt}\lineto[r]&\undertag{\bullet}{P'(-2)}{4pt}
&&&&\\
&&&&\\
&&&&\\
&&&&
}
$$
Again, one can blow-up between $Q'$ and $N'$ repeatedly, and have a completion of $E'$ by adjoining a $-1$ at the end of $\mathcal C(P)$ with the $-(p+2)$ and the end of $\mathcal C(Q)$ with $Q$.  If however $p=1$, one may blow up between $Q'$ and $M'$, and find another solution by letting $-1$'s attach $\mathcal C(P)$ with the $-2$ curve to the right of $Q$, and $\mathcal C(Q)$ with $Q$.   This gives the desired existence of only one (or two, in the special case) configuration(s) of $-1$'s on $Z'$, and then pulling back  gives Basic Models on $Z$.

\subsection{Type $\Gamma=\mathcal M(0,q,0)$}

One uses the same $Z'$.  There is now one chain, $\mathcal C(Q)$, which consists solely of $-2$ curves.  One makes it blow-downable only by letting a $-1$ connect $Q$ itself with one or the other end of the chain.  It is easy to check that this is possible only if the $-1$ is at the far end of the chain, so there is a unique way to blow-down.  So, one must be in the case of the Basic Model.

\section{Self-isotropic subgroups and Fowler's method}

If $L$ is a non-degenerate lattice, the dual $L^*=\text{Hom}(L,\mathbb Z)$ admits a non-degenerate pairing into $\mathbb Q$.  Thus, the finite  \emph{discriminant group} $D(L)=L^*/L$ admits a  non-degenerate \emph{discriminant pairing} into $\mathbb Q/\mathbb Z$.   Overlattices $L\subset M$ of the same rank correspond to isotropic subgroups $\bar{M}=M/L$ of the discriminant group.  If $M$ is unimodular, then $\bar{M}$ is self-isotropic.  The importance of these notions in smoothing surface singularities may be found in \cite{lw}, Section 2.

 If $\Gamma$ is one of the dual graphs listed in the Bhupal-Stipsicz Theorem, it gives rise to a lattice and discriminant group $D(\Gamma)$.  On a $\Gamma$ surface $(Z,E)$, the lattice $\mathbb E(\Gamma)=\oplus_i \mathbb Z [E_i]$ spanned by the divisor classes in Pic($Z$) comes with an identification with the lattice of $\Gamma$.  
 
 \begin{definition}  Let $(Z,E)$ be a $\Gamma$ surface.  Then the \emph{self-isotropic group of $(Z,E)$} is the subgroup $I$ of $D(\Gamma)$ associated with the unimodular overlattice $\mathbb E(\Gamma)\subset Pic(Z)\subset  \mathbb E(\Gamma)^*$.   
 \end{definition}

  Fowler studied in \cite{jf} the important map
  $$\xi:\{\text{Isomorphism classes of}\   \Gamma\  \text{surfaces} \ (Z,E)\}\rightarrow \{\text{Self-isotropic subgroups of}\  D(\Gamma)\}$$

\subsection{Fowler's approach}

\begin{definition}  A self-isotropic subgroup $I\subset D(\Gamma)$ is called \emph{basic} if it is associated with a Basic Model $\Gamma$ surface. 
\end{definition}

Whenever $D(\Gamma)$ has only one self-isotropic subgroup (as happens most of the time), then of course being ``basic'' is not an extra condition. On the other hand, there are examples of  $\Gamma$ of types $\mathcal W, \mathcal N, \mathcal M$ which have non-basic self-isotropic subgroups (hence the need for the Theorem in this paper). 
One of Fowler's main theorems is

\begin{theorem}\cite{jf}  Let $(Z,E)$ be a $\Gamma$ surface whose self-isotropic subgroup is basic.  Then $(Z,E)$ is basic.
\end{theorem}

If a $\Gamma$ surface $(Z,E)$ has basic self-isotropic subgroup $I$, the goal is to prove that it is itself basic.  This is accomplished by locating in a precise location $-1$ curves on $Z$ which allow one to blow down in a unique way.  The assumption that $I$ is basic implies that for each potential $-1$ curve, there is a line bundle $\mathcal L$ with the correct intersection properties with all $E_i$.  If $\mathcal L$ has a section giving an irreducible curve, it will be the sought after $-1$.  


 Fowler achieves this for every potential $-1$ curve, via a case by case look at all the types of $\Gamma$, of course assuming that the self-isotropic subgroup is basic.  

Fowler starts with general considerations about $K$ similar to those in Section $1$ above.  He concludes (using Riemann-Roch) that $\mathcal L$ has a non-zero section, so is represented by an effective divisor $$L=\sum n_iE_i + \sum m_jF_j.$$
Here the $F_j$ are irreducible curves not among the $E_i$.  One needs to show that all the $n_i$ are $0$ as well as all but one of the $m_j$.  

The method involves locating in each case various exceptional nef divisors $N$.  Then $$N\cdot L=\sum n_i(N\cdot E_i)+\sum m_j(N\cdot F_j)\geq 0.$$ By $L$'s intersection properties, one can frequently arrange that the product  $N\cdot L=0$.  In this case, if $N\cdot E_i>0$, then necessarily $n_i=0$.  Note that if $F_j$ is not a $-1$ curve, then it is nef (Corollary 1.3).

For instance, in case $\mathcal W$, suppose one wishes to prove the existence of a $-1$ curve connecting the end of the chain $\mathcal C(P)$ with $R$.  Let $L$ be the divisor above representing the potential curve, and choose first $N=E_0$, the central curve.  Then $L\cdot N=0$, while $N$ dots to $1$ with the central curve and the three adjacent curves.  We conclude that the corresponding four coefficients $n_i$ in the expansion of $L$ equal $0$.  If we choose as nef divisor $N=(p+2)E_0+P$, then again $L\cdot N=0$, so the coefficient of the neighbor of $P$ is $0$ as well.  

Fowler develops very efficient methods for all $\Gamma$ to show that each potential $-1$ curve actually does exist.  This allows careful analysis of the blow-down.  

\subsection{Number of $\Q$HD smoothing components}
Combining the main result Theorem $0.6$ of this paper with Fowler's results, here is the final count of $\Q$HD smoothing components for weighted homogeneous singularities:
\begin{enumerate}
\item  Two components for $\mathcal W(p,p,p)$, $\mathcal N(q+2, q,0)$, and $\mathcal M(r+1,q,r)$, with two different self-isotropic subgroups in each case.
\item A unique component for all other $\mathcal W, \mathcal N$, and $\mathcal M$.
\item A unique component for type $\mathcal B^3_2$ and $\mathcal C^3_3$.
\item Two components for type $\mathcal C^3_2$, with the same self-isotropic subgroup.
\item Two components for type $\mathcal A^4$ with two different isotropic subgroups in each case.
\item Two components for types $\mathcal B^4$ and $\mathcal C^4(p), p>0$, with one self-isotropic subgroup in each case.
\item A unique component for type $\mathcal C^4(0)$.
\end{enumerate}

As mentioned before,  Fowler shows the existence of two components is a consequence either of a symmetry in the graph $\Gamma$ or of complex conjugation in the blowing-up process.

\includepdf[pages=1, pagecommand={}, offset=-0.3cm -2cm]{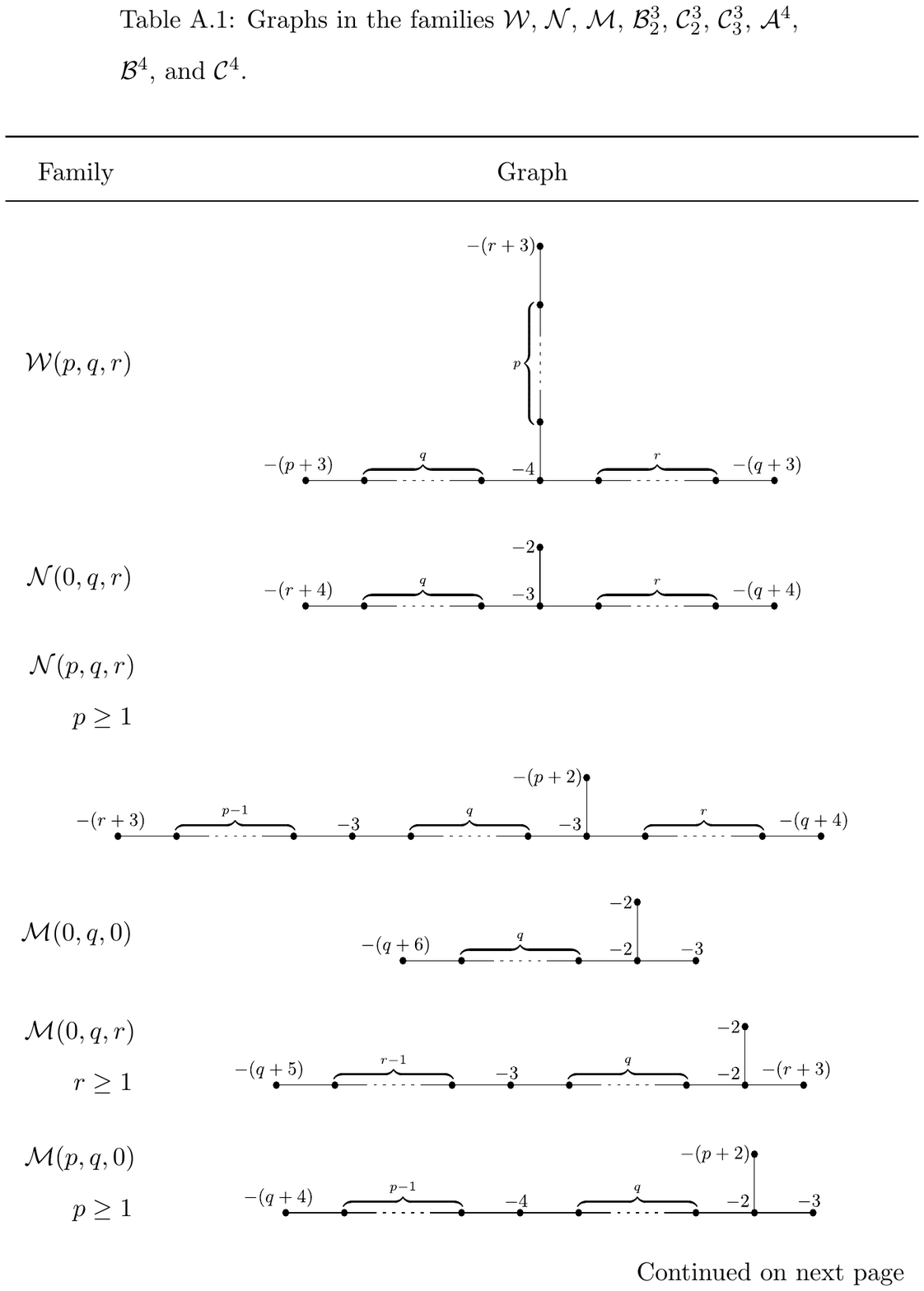}
\includepdf[pages=2, pagecommand={}, offset=0.3cm -2cm]{Fowlerphd.pdf}
\includepdf[pages=3, pagecommand={}, offset=-0.3cm -2cm]{Fowlerphd.pdf}
\includepdf[pages=4, pagecommand={}, offset=0.3cm -2cm]{Fowlerphd.pdf}
\includepdf[pages=5, pagecommand={}, offset=-0.3cm -5 cm]{Fowlerphd.pdf}


\end{document}